\documentclass[12pt,reqno,letter]{amsart}

\usepackage[all]{xy}
\usepackage{texdraw}

\def\theoremnumberstyle{plain}
\def\proofname{Proof}
\usepackage{tom}

\newcommand{\tom}[1]{}

\newcommand{\pp}{{\mathbb P}}

\SelectTips{cm}{12}

\begin{document}

   \parindent0cm

   \title[Triple Point]{Triple-Point Defective Ruled Surfaces}
   \author{Luca Chiantini}
   \address{Universit\'a degli Studi di Siena\\
     Dipartimento di Matematica\\
     Via del Capitano, 15\\
     I -- 53100 Siena
     }
   \email{chiantini@unisi.it}
   \urladdr{http://pcserver.mat.unisi.it/web/chiantini/chlam.htm}
   \author{Thomas Markwig}
   \address{Universit\"at Kaiserslautern\\
     Fachbereich Mathematik\\
     Erwin-Schr\"odinger-Stra\ss e\\
     D -- 67663 Kaiserslautern
     }
   \email{keilen@mathematik.uni-kl.de}
   \urladdr{http://www.mathematik.uni-kl.de/\textasciitilde keilen}
   \thanks{The second author was supported by the EAGER node of Torino
     and by the Institute for Mathematics and its Applications (IMA)
     in Minneapolis.} 


   \date{December, 2006.}


   \begin{abstract}
     In \cite{CM07} we studied triple-point defective very ample linear systems
     on regular surfaces, and we showed that they can only exist if the
     surface is ruled. In the present paper we show that we can drop
     the regularity assumption, and we classify the
     triple-point defective very ample linear systems on
     ruled surfaces.
   \end{abstract}

   \maketitle

   Let $S$ be a smooth projective surface,
   $K=K_S$  the canonical class and $L$ a
   divisor class on $S$ 

   We study a classical interpolation problem for the pair $(S,L)$,
   namely whether for a general point $p\in S$ the linear system
   $|L-3p|$ has the expected dimension
   \begin{displaymath}
     \expdim|L-3p|=\max\{-1,\dim|L|-6\}.
   \end{displaymath}
   If this is not the case we call the pair $(S,L)$
   \emph{triple-point defective}. 

   This paper is indeed a continuation of \cite{CM07}, where some classification
   of triple point defective pairs is achieved, under the assumptions:

   \begin{displaymath}
     L, L-K \mbox{ \it very ample, and }  (L-K)^2>16, 
   \end{displaymath}
   conditions that we will take all over the paper.

   With these assumptions, the main result of \cite{CM07} says that
   all triple-point defective {\it regular} surfaces are rationally ruled.

   We tackled the problem by considering $|L-3p|$ as fibres of the the map
   $\alpha$ in the following diagram,
   \begin{equation}\label{eq:alphabeta}
     \xymatrix{
       |L|=\pp(H^0(L)^*)&\kl_3\ar[r]^\alpha\ar[l]_(.25)\beta & S       
     }
   \end{equation}
   where $\kl_3$ denotes the incidence variety
   \begin{displaymath}
     \kl_3=\{(C,p)\in|L|\times S\;|\;\mult_{p}(C)\geq 3\}
   \end{displaymath}
   and $\alpha$ and $\beta$ are the obvious projections.

   Assuming that for a general point $p\in S$ there is a curve in
   $L_p$ with a triple-point in $p$ -- and hence $\alpha$ surjective,
   we considered then the {\it equimultiplicity scheme} 
   $Z_p$ of a curve $L_p\in|L-3p|$ defined by 
   \begin{displaymath}
     \kj_{Z_p,p}=\left\langle\frac{\partial
         f_p}{\partial x_p},\frac{\partial f_p}{\partial
         y_p}\right\rangle + \langle x_p,y_p\rangle^3.
   \end{displaymath}
   One easily sees that $(S,L)$ triple-point defective necessarily
   implies that 
   \begin{displaymath}
     h^1\big({S},\kj_{Z_p}(L)\big)\not=0.
   \end{displaymath}
   
   Non--zero elements in $H^1\big({S},\kj_{Z_p}(L)\big)$ determine by Serre duality
   a non--trivial extension $\ke_p$ of $\kj_{Z_p}(L-K)$ by  $\ko_S$, which 
   turns out to be a rank $2$ bundle on the surface. Due to
   the assumption $(L-K)^2>16$,  $\ke_p$ is Bogomolov unstable. We then
   exploited the destabilizing divisor $A_p$ of $\ke_p$ in order to
   obtain the above mentioned result. 
   \smallskip
   
   For non--regular surfaces, the argument of \cite{CM07} shows the following
   lemma (see \cite{CM07}, Prop.\ 17 and Prop.\ 18):
   
   \begin{proposition}\label{prop:length4}
     Suppose that, with the notation
     in \eqref{eq:alphabeta}, $\alpha$ is surjective,
     and suppose as usual that $L$ and $L-K$ are
     very ample with $(L-K)^2>16$. 
     
     For $p$ general in $S$ and for $L_p\in |L-3p|$ general, call
      $Z_p'$  the minimal subscheme of the equimultiplicity
     scheme $Z_p$ of $L_p$ such that  $$h^1\big(S,\kj_{Z'_p}(L)\big)\not=0.$$
     Then either:
     \begin{enumerate}
     \item[1)] $\length(Z_p')=3$ and $S$ is ruled; or
     \item[2)] $\length(Z_p')=4$ and, for $p\in S$ general, 
     there are smooth, elliptic curves
     $E_p$ and $F_p$ in $S$ through $p$ such that $E_p^2=F_p^2=0$,
     $E_p.F_p=1$ and $L.E_p=L.F_p=3$.
     In particular, both $|E|_a$ and $|F|_a$  induce an elliptic
     fibration with section on $S$ over an elliptic  curve.
     \end{enumerate}
   \end{proposition}
   
   This is our starting point. We will in this paper show that the latter case
   actually cannot occur, and we will classify the triple-point
   defective linear systems $L$ as above on  ruled
   surfaces. It will in particular follow that the fibre of the
   ruling is contained exactly twice, and thus that the map $\beta$
   above is generically finite. 
   
   Our main results are:
   
   \begin{theorem}\label{thm:main}
     Suppose that the pair $(S,L)$ is triple-point defective where
     $L$ and $L-K$ are very ample with $(L-K)^2>16$.
     Then $S$ admids a ruling $\pi:S\rightarrow C$.
   \end{theorem}

   For the classification, call $C_0$ a section of the ruled surface $S$,
   $\mathfrak{e}$ the line bundle on the base curve given by
   the determinant of the defining bundle, and call $E_i$ the exceptional divisors
   (see pp.\ \pageref{sec:ruled} and \pageref{page:notation} for a
   more precise setting of the notation):  
   
   \begin{theorem}\label{main2} 
     Assume that $\pi:S\rightarrow C$ is a ruled surface
     and that the pair $(S,L)$ is triple-point defective, where
     $L$ and $L-K$ are very ample with $(L-K)^2>16$. 

     Then $\pi$ is minimal, i.e.\ $S$ is geometrically ruled, and  for
     a general point $p\in S$ the linear system $|L-3p|$ contains  
     a fibre of the ruling as fixed component with  multiplicity two. 

     Moreover, in the previous notation, the line bundle 
     $L$ is of type $C_0+\pi^*\mathfrak{b}$ for some
     divisor $\mathfrak{b}$ on $C$ such that
     $\mathfrak{b}+\mathfrak{e}$ is very ample.  
   \end{theorem}

   In Section \ref{sec:product} we
   will first show that a surface $S$ admitting two elliptic
   fibrations as required by Proposition \ref{prop:length4} would necessarily
   be a product of two elliptic curves and the triple-point defective
   linear system would be of type $(3,3)$. We then show that such a
   system is never triple-point defective, setting the first part of the main
   theorem. 
   
   In Section \ref{sec:ruled}
   we  classify the triple-point defective linear systems on
   ruled surfaces, thus completing our main results.


   \section{Products of Elliptic Curves}\label{sec:product}

   In the above setting, consider a triple-point
   defective tuple $(S,L)$ where the equimultiplicity
   scheme $Z_p$ (see \cite{CM07}) of a general element $L_p\in|L-3p|$
   \  admitted a
   complete intersection subscheme $Z_p'$ of length \emph{four} with 
   \begin{displaymath}
     h^1\big(S,\kj_{Z'_p}(L)\big)\not=0.
   \end{displaymath} 
   As explained in the introduction, Prop. \ref{prop:length4}, after
   \cite{CM07} we      know that, for $p\in S$ general, there are
   smooth, elliptic curves 
   $E_p$ and $F_p$ in $S$ through $p$ such that $E_p^2=F_p^2=0$,
   $E_p.F_p=1$ and $L.E_p=L.F_p=3$.

   In particular, both $|E|_a$ and $|F|_a$  induce an elliptic
   fibration with section on $S$ over an elliptic  curve.

   We will now show that this situation indeed cannot occur.
   Namely, for general $p$ and $L_p$ there cannot
   exist such a scheme $Z_p'$. 

   \begin{lemma}\label{lem:ellipticfibrations}
     Suppose that the surface $S$ has two elliptic fibrations
     $\pi:S\longrightarrow E_0$ and $\pi':S\longrightarrow F_0$ with
     general fibre $E$ respectively $F$ satisfying
     $E.F=1$. 

     Then $E_0$ and $F_0$ are elliptic curves, and $S$ is the blow-up
     of a product of two elliptic curves $S'=E\times E_0\cong E\times F$. 

   \end{lemma}
   \begin{proof}
     Since $E.F=1$ we have that $F$ is a section of $\pi$, and
     thus $F\cong E_0$ via $\pi$. In particular, $E_0$ and, similarly,
     $F_0$ are elliptic curves.

     It is well known that there are no non--constant maps from a rational
     curve to a curve of positive genus (\cite{Har77}, IV.2.5.4). 
     Thus any exceptional curve of $S$ 
     sits in some fiber. Thus we can reach relatively
     minimal models of $\pi$   and $\pi'$ by successively blowing down
     exceptional $-1$-curves which  belong to fibres of both $\pi$ and
     $\pi'$, i.e.\ we have
     the following commutative diagram
     \begin{displaymath}
       \xymatrix{
         S\ar[dr]^(.7){\phi}\ar@/^/[drr]^{\pi}\ar@/_/[ddr]_{\pi'}& \\
         & S'\ar[r]^(.4){\widetilde{\pi}}\ar[d]^{\widetilde{\pi}'} &E_0\\
         & F_0 &
       }
     \end{displaymath}
     where $S'$ is actually a minimal surface. Since a general fibre
     of $\pi$ or $\pi'$ is not touched by the blowing-down $\phi$ we
     may denote the general fibres of $\widetilde{\pi}$ and $\widetilde{\pi}'$
     again by $E$ respectively $F$, and we still have $E.F=1$.

     We will now try to identify the minimal surface $S'$ in the
     classification of minimal surfaces.

     By \cite{Fri98} Ex.\ 7.9 the canonical divisor $K_{S'}$ is
     numerically trivial, since $S'$ is a minimal surface admitting
     two elliptic fibrations over elliptic curves. 

     But then we can apply \cite{Fri98} Ex.\ 7.7, and since the base curve
     $E_0$ of the fibration $\widetilde{\pi}$ is
     elliptic we see  that the invariant
     $d=\deg(L)=\deg\big((R^1\pi_*\ko_{S'})^{-1}\big)$ of the relatively minimal
     fibration $\widetilde{\pi}$ mentioned in \cite{Fri98} Cor.\ 7.17 is 
     zero, so that the  same corollary implies that the fibration
     has at most multiple fibres with smooth reduction as singular
     fibres. However, since $\widetilde{\pi}$ has a section $F$
     there are no multiple fibres, and thus all fibres of
     $\widetilde{\pi}$ are smooth. 

     Moreover, since the canonical divisor of $S'$ is numerically
     trivial it is in particular nef, and by \cite{Fri98} Thm.\ 10.5 we
     get that the Kodaira dimension $\kappa(S')$ of $S'$ is zero. 

     Moreover, by \cite{Fri98} Cor.\ 7.16 the surface $S'$ has second
     Chern class $c_2(S')=0$, since the invariant
     $d=\deg\big((R^1\pi_*\ko_{S'})^{-1}\big)=0$ as already mentioned
     above. Thus by the Enriques-Kodaira Classification (see e.g.\
     \cite{BHPV04} Thm.\ 10.1.1) $S'$ must either be a torus or
     hyperelliptic (where the latter is sometimes also called
     bielliptic). A bielliptic surface has precisely two elliptic
     fibrations, but one of them is a fibration over a $\pp^1$ and
     only one is over an elliptic curve (see e.g.\ \cite{Rei97}
     Thm.\ E.7.2). Thus $S'$ is not bielliptic. Moreover, if $S'$ is
     a torus then $K_{S'}$ is  
     trivial and thus so is $(R^1\pi_*\ko_{S'})^{-1}$, which by
     \cite{Fri98} Cor.\ 7.21 implies that $S'$ is a product of the
     base curve with a fibre.
   \end{proof}

   Lemma \ref{lem:ellipticfibrations} implies that in order to show that
   the situation of Proposition \ref{prop:length4} cannot occur, we
   have to understand products of elliptic curves. 

   Let us, therefore, consider a surface ${S}=C_1\times C_2$ which is
   the product of two smooth elliptic curves.
   
   Let us set some notation. We will use some results by \cite{Kei01} 
   Appendices G.b and G.c in the sequel.
    
   The surface ${S}$ is naturally equipped with two projections
   $\pi_i:{S}\longrightarrow C_i$. If $\mathfrak{a}$ is a divisor on
   $C_2$ of degree $a$ and $\mathfrak{b}$ is a divisor on $C_1$ of
   degree $b$ then the divisor
   $\pi_2^*\mathfrak{a}+\pi_1^*\mathfrak{b}\sim_a aC_1+bC_2$, where by
   abuse of notation we denote by $C_1$ a fixed fibre of $\pi_2$ and
   by $C_2$ a fixed fibre of $\pi_1$. Moreover, $K_{S}$ is trivial, and
   given two divisors $D\sim_a aC_1+bC_2$ 
   and $D'\sim_a a'C_1+b'C_2$ then the intersection product is
   \begin{displaymath}
     D.D'=(aC_1+bC_2).(a'C_1+b'C_2)=a\cdot b'+a'\cdot b.
   \end{displaymath}

We will consider first the case 
$$L=\pi_2^*(\mathfrak{a})+\pi_1^*(\mathfrak{b})$$
 where both $\mathfrak{b}$ on $C_1$ and $\mathfrak{a}$ on $C_2$
are divisors of degree $3$. The dimension of the linear system $|L|$ is 
$ \dim|L|=8,$  and thus for a point $p\in {S}$ the expected dimension
is $\expdim|L-3p|=\dim|L|-6=2$.

Notice that a divisor of degree three on an elliptic curve is always
   very ample and embeds the curve as a smooth cubic in $\pp^2$. Since
   the smooth plane cubics are classified by their normal forms
   $xz^2-y\cdot(y-x)\cdot(y-\lambda\cdot x)$ with $\lambda\not=0$ the
   following example reflects the behaviour of any product of elliptic
   curves embedded via a linear system of bidegree $(3,3)$.
   
 \begin{example}\label{ex:ellipticcomputation}
     Consider two smooth plane cubics 
     \begin{displaymath}
       C_1=V\big(xz^2-y\cdot (y-z)\cdot(y-az)\big)
     \end{displaymath}
     and
     \begin{displaymath}
       C_2=V\big(xz^2-y\cdot (y-z)\cdot(y-bz)\big).
     \end{displaymath}
     The surface ${S}=C_1\times C_2$ is embedded into $\pp^8$ via the
     Segre embedding
     \begin{displaymath}
       \phi:\pp^2\times\pp^2\longrightarrow\pp^8:
       ((x_0:x_1:x_2),(y_0:y_1:y_2))
       \mapsto
       (x_0y_0:\ldots:x_2y_2).
     \end{displaymath}
     We may assume that both curves contain the point $p=(1:0:0)$ as a general 
     non-inflexion point, and the point $(p,p)$ is mapped by the Segre embedding to
     $\phi(p,p)=(1:0:\ldots:0)$.  If we denote by $z_{i,j}$,
     $i,j\in\{0,1,2\}$, the coordinates on $\pp^8$ as usual, then
     the maximal ideal locally at $\phi(p,p)$ is generated by $z_{0,2}$ and
     $z_{2,0}$, i.e.\ these are local coordinates of ${S}$ at
     $\phi(p,p)$. A standard basis computation shows that locally at
     $\phi(p,p)$ the coordinates $z_{i,j}$ satisfy modulo the ideal
     of ${S}$ and up to  multiplication by a unit the following
     congruences (note, $z_{0,0}=1$)
     \begin{align*}
       z_{0,1}\equiv & \frac{1}{b}\cdot z_{0,2}^2,&
       z_{1,0}\equiv & \frac{1}{a}\cdot z_{2,0}^2,&
       z_{1,1}\equiv & \frac{1}{ab}\cdot z_{0,2}^2\cdot z_{2,0}^2,\\
       z_{1,2}\equiv & \frac{1}{a}\cdot z_{0,2}\cdot z_{2,0}^2,&
       z_{2,1}\equiv & \frac{1}{b}\cdot z_{0,2}^2\cdot z_{2,0},&
       z_{2,2}\equiv & z_{0,2}\cdot z_{2,0}.
     \end{align*}
     Thus a
     hyperplane section $H=a_{0,0}z_{0,0}+\ldots+a_{2,2}z_{2,2}$
     of ${S}$ is locally in $\phi(p,p)$ modulo 
     $\m^3=\langle z_{0,2},z_{2,0}\rangle^3$ given by
     \begin{displaymath}
       H\equiv a_{0,0}+a_{0,2}z_{0,2}+a_{2,0}z_{2,0}
       +\frac{a_{0,1}}{b}\cdot z_{0,2}^2
       +\frac{a_{1,0}}{a}\cdot z_{2,0}^2
       +a_{2,2}z_{0,2}z_{2,0},
     \end{displaymath}
     and hence the family of hyperplane sections having multiplicity
     at least three in $\phi(p,p)$ is given by 
     \begin{displaymath}
       a_{0,0}=a_{0,1}=a_{1,0}=a_{0,2}=a_{2,0}=a_{2,2}=0.
     \end{displaymath}
     But then the family has parameters $a_{1,1},a_{1,2},a_{2,1}$, and its
     dimension coincides with the expected dimension $2$. Moreover,
     the $3$-jet of a hyperplane section $H$ through $\phi(p,p)$ with
     multiplicity at least three is
     \begin{displaymath}
       \jet_3(H)\equiv z_{0,2}\cdot
       z_{2,0}\cdot\left(\frac{a_{1,2}}{a}\cdot
         z_{2,0}+\frac{a_{2,1}}{b}\cdot z_{0,2}\right),
     \end{displaymath}
     which shows that for a general choice of $a_{2,1}$ and $a_{1,2}$
     the point $\phi(p,p)$ is an ordinary triple point.
   \end{example}

   \begin{remark}\label{ex:elliptic}
     We actually can say very precisely what it means that $p$ is
     general in the product, namely that neither $\pi_1(p)$ is a inflexion point of
     $C_1$, nor $\pi_2(p)$ is a inflexion point of $C_2$. 
     
     Indeed, since $\mathfrak{a}$ is very ample of degree three, for each
     point $p\in {S}$ there is a unique point  $q_a\in C_2$ such that 
     $q_a+2\cdot\pi_2(p)\sim_l\mathfrak{a}$.
     When $\pi_2(p)$ is a inflexion point of $C_2$, then
     $q_a=\pi_2(p)$ and thus the two-dimensional family
     \begin{displaymath}
       3C_{1,\pi_2(p)}+|\pi^*(\mathfrak{b})|\subset|L-3p|
     \end{displaymath}
     gives a superabundance of the dimension of $|L-3p|$ by one.
     
     Similarly one can argue when $\pi_1(p)$ is a inflexion point of $C_1$.
   \end{remark}
      
   
   Now we are ready for the proof of Theorem \ref{thm:main}.
   
   \begin{proof}[Proof of Theorem \ref{thm:main}]
     By Proposition \ref{prop:length4}, it is enough to prove that
     when $S$ has two elliptic fibrations as in the proposition, then
     $S$ is not triple--point defective. 
  
     By Lemma \ref{lem:ellipticfibrations}, $S$ is the blow-up
     $\pi:S\longrightarrow S'$ of a
     product $S'=C_1\times C_2$ of two elliptic curves, and we may
     assume that the curves $E_p$ and $F_p$ in Proposition
     \ref{prop:length4} are the fibres of $\pi_1$ respectively
     $\pi_2$.

     Our first aim will be to show that actually $S=S'$. For this note
     that 
     \begin{displaymath}
       \Pic(S)=\bigoplus_{i=1}^k E_i\oplus\pi^*\Pic(S'),
     \end{displaymath}
     where the $E_i$ are the total transforms of the exceptional
     curves arising throughout the blow-up, i.e.\ the $E_i$ are (not
     necessarily irreducible) rational curves with self-intersection
     $E_i^2=-1$ and such that $E_i.E_j=0$ for $i\not=j$ and
     $E_i.\pi^*(C)=0$ for any curve $C$ on $S'$. In particular, since
     $K_{S'}$ is trivial we have that $K_S=\sum_{i=1}^k E_i$, and if
     $L=\pi^*L'-\sum_{i=1}^ke_iE_i$ then
     $L-K=\pi^*L'-\sum_{i=1}^k(e_i+1)E_i$. We therefore have
     \begin{displaymath}
       16<(L-K)^2=(L')^2-\sum_{i=1}^k (e_i+1)^2,
     \end{displaymath}
     or equivalently
     \begin{equation}\label{eq:corlength4:1}
       (L')^2\geq 17+\sum_{i=1}^k(e_i+1)^2\geq 17+4k,
     \end{equation}
     where the latter inequality is due to the fact that $e_i=L.E_i>0$
     since $L$ is very ample.
     By the assumption of Proposition \ref{prop:length4} we know that
     $L'.C_1=L.E_p=3$ and $L'.C_2=L.F_p=3$, and therefore by
     \cite{Har77} Ex.\ V.1.9 
     \begin{equation}\label{eq:corlength4:2}
       (L')^2\leq 2\cdot (L'.C_1)\cdot (L'.C_2)=18.
     \end{equation}
     But \eqref{eq:corlength4:1} and \eqref{eq:corlength4:2} together
     imply that no exceptional curve exists, i.e.\ $S=S'$.

     Since now $S$ is a product of two elliptic curves, by \cite{LB92}
     we know that the Picard number $\rho=\rho(S)$ satisfies $2\leq\rho\leq 4$, and
     the N\'eron-Severi group can be generated by the two general
     fibres $C_1$ and $C_2$ together with certain graphs $C_j$, $3\leq
     j\leq \rho$, of morphisms $\varphi_j:C_1\longrightarrow C_2$. In particular,
     $C_j.C_2=1$ and $C_j.C_1=\deg(\varphi_j)\geq 1$ for $3\leq j\leq
     \rho$. Moreover, these graphs have self intersecting zero. If we
     now assume that $L\sim_a \sum_{j=1}^\rho a_iC_i$ then
     \begin{displaymath}
       L^2=2\cdot \sum_{i<j} a_i\cdot a_j\cdot (C_i.C_j)
     \end{displaymath}
     is divisible by $2$, and since $L=L-K$ with $(L-K)^2>16$ we
     deduce with \cite{Har77} Ex.\ V.1.9 that
     \begin{displaymath}
       L^2=(L-K)^2=18=2\cdot (L.C_1)\cdot (L.C_2),
     \end{displaymath}
     and thus that
     \begin{displaymath}
       L\sim_a 3C_1+3C_2,
     \end{displaymath}
     or in equivalently, that 
     \begin{displaymath}
       L=\pi_2^*\mathfrak{a}+\pi_1^*\mathfrak{b}
     \end{displaymath}
     for some divisors $\mathfrak{a}$ on $C_2$ and $\mathfrak{b}$ on
     $C_1$, both of degree $3$. That is, we are in the situation of
     Example \ref{ex:ellipticcomputation}, and we showed there that $(S,L)$
     then is not triple-point defective.  
   \end{proof}

   \begin{remark}\label{remlength} Notice that, in practice,  since
     \begin{displaymath}
       h^1(S,L)=h^0(C_1,\mathfrak{b})\cdot h^1(C_2,\mathfrak{a})
       +h^0(C_2,\mathfrak{a})\cdot h^1(C_1,\mathfrak{b})=0,
     \end{displaymath}
     the non-triple-point defectiveness  shows that for general 
     $p\in S$ and $L_p\in |L-3p|$ no 
     $Z'_p$ as in the assumptions of Proposition \ref{prop:length4}
     can have length $4$.
   \end{remark}


   \section{Geometrically Ruled Surfaces}\label{sec:ruled}

   Let $\xymatrix@C0.6cm{S=\pp(\ke)\ar[r]^(0.65)\pi & C}$ be a
   geometrically ruled surface with normalized bundle $\ke$ (in the
   sense of \cite{Har77} V.2.8.1). The N\'eron-Severi group of $S$ is 
   \begin{displaymath}
     \NS(S) = C_0\Z\oplus f\Z,
   \end{displaymath}
   with intersection matrix
   \begin{displaymath}
     \left(\!\begin{array}{rc}-e & 1 \\ 1 & 0\end{array}\right),
   \end{displaymath}
   where $f\cong\pp^1$ is a fixed fibre of $\pi$, $C_0$ a fixed section of $\pi$
   with $\ko_S(C_0)\cong\ko_{\pp(\ke)}(1)$,
    and
   $e=-\deg(\mathfrak{e})\geq -g$ where
   $\mathfrak{e}=\Lambda^2\ke$.\tom{\footnote{By \cite{Nag70} Theorem 1 
     there is some section $D\sim_a C_0+bf$ with $g\geq
     D^2=2b-e$. Since $D$ is irreducible, by \cite{Har77} V.2.20/21
     $b\geq 0$, and thus $-g\leq e$.}} 
   If
   $\mathfrak{b}$ is a divisor on $C$ we will write $\mathfrak{b}f$
   for the divisor $\pi^*(\mathfrak{b})$ on $S$, and so for the
   canonical divisor we have 
   \begin{displaymath}
     K_S \;\;\sim_l\;\; -2C_0+(K_C+\mathfrak{e})\cdot f\;\;\sim_a\;\; -2C_0+ (2g-2-e)f,
   \end{displaymath}
   where $g=g(C)$ is the genus of the base curve $C$.

   \begin{example}\label{ex:ruled}
     Let $\mathfrak{b}$ be a divisor on $C$ such that $\mathfrak{b}$
     and $\mathfrak{b+e}$ are both very ample and such that
     $\mathfrak{b}$ is non-special. If $C$ is rational we should in
     addition assume that $\deg(\mathfrak{b})+\deg(\mathfrak{b+e})\geq
     6$.
     Then the divisor
     $L=C_0+\mathfrak{b}f$ is very ample (see e.g.\ \cite{FuP00} Prop.\
     2.15) of dimension
     \begin{displaymath}
       \dim|L|=h^0(C,\mathfrak{b})+h^0(C,\mathfrak{b}+\mathfrak{e})-1
     \end{displaymath}
     Moreover, for any point $p\in S$ we then have (see \cite{FuP00}
     Cor.\ 2.13) 
     \begin{displaymath}
       \dim|\, C_0+(\mathfrak{b}-2\pi(p))\cdot f|=\dim|\, C_0+\mathfrak{b}f|-4
       =h^0(C,\mathfrak{b})+h^0(C,\mathfrak{b}+\mathfrak{e})-5,
     \end{displaymath}
     and we have for $p$ general
     \begin{displaymath}
       \dim|\, C_0+(\mathfrak{b}-2\pi(p))\cdot f-p\, |=
       h^0(C,\mathfrak{b})+h^0(C,\mathfrak{b}+\mathfrak{e})-6.
     \end{displaymath}
     For this note that $\mathfrak{b}$ and $\mathfrak{b+e}$ very ample
     implies that this number is non-negative -- in the rational case
     we need the above degree bound.
     
     If we denote by $f_p=\pi^*\big(\pi(p)\big)$ the fibre of $\pi$
     over $\pi(p)$, then by B\'ezout and since $L.f_p=(L-f_p).f_p=1$ we
     see that $2f_p$ is a fixed component of $|L-3p|$ and we have
     \begin{displaymath}
       |L-3p|=2f_p+|\, C_0+(\mathfrak{b}-2\pi(p))\cdot f-p\, |,
     \end{displaymath}
     so that
     \begin{multline*}
       \dim|L-3p|=h^0(C,\mathfrak{b})+h^0(C,\mathfrak{b}+\mathfrak{e})-6
       =\dim|L|-5\\>\dim|L|-6=\expdim|L-3p|.
     \end{multline*}
     This shows that $(S,L)$ is triple-point defective and $|L-3p|$
     contains a fibre of the ruling as double component. Moreover,
     for a general $p$ the linear series $|L-3p|$ cannot contain a
     fibre of the ruling more than twice due to the above dimension
     count for $|\, C_0+(\mathfrak{b}-2\pi(p))\cdot f-p\, |$.
   \end{example}
     
   Next we are showing that a geometrically ruled surface is indeed
   triple-point defective with respect to a line bundle $L$ which
   fulfills our assumptions, and in Corollary
   \ref{cor:geometricallyruled} we will see 
   that this is not the case for non-geometrically ruled surfaces.
   
   \begin{proposition}  
     On every geometrically ruled surface 
     $S=\pp(\ke)\stackrel{\pi}{\longrightarrow} C$
     there exists some very ample line bundle $L$ such that the pair 
     $(S,L)$ is triple--point defective, and moreover also $L-K$ is
     very ample with $(L-K)^2>16$.  
   \end{proposition}
   \begin{proof} 
     It is enough to take
     $L=C_0+\mathfrak{b}f$, with $b=\deg(\mathfrak{b})=3a$ such that
     $a, a-e, a+e, a-2g+2+e, a-2g+2-e$ are all bigger or equal than $2g+1$.
     
     Indeed in this case $\mathfrak{b}$ and  $\mathfrak{b}+\mathfrak{e}$ 
     are both very ample. For $p\in C$ general, we also have
     that both $\mathfrak{b}-p$ and $\mathfrak{b}+\mathfrak{e}-p$ are
     non-special. It follows that $L$ is very ample (by \cite{Har77} Ex.\ V.2.11.b)
     and $(S,L)$ is triple point defective, by the previous example.
     Moreover, in this situation we have:
     \begin{displaymath}
       L-K\sim_l 3C_0+\big(\mathfrak{b}-K_C-\mathfrak{e}\big)\cdot f.
     \end{displaymath}
     Hence      
     \begin{displaymath}
       (L-K)^2=\big(3C_0+(\deg(\mathfrak{b})-2g+2+e)\cdot f\big)^2
       \geq 18>16.
     \end{displaymath}
     Finally, if we fix a divisor $\mathfrak{a}$ of degree $a$ on $C$,
     then $L-K$ is the sum of the divisors
     $C_0+\big(\mathfrak{a}-K_C\big)\cdot f$,
     $C_0+\big(\mathfrak{a}-\mathfrak{e}\big)\cdot f$,
     $C_0+\mathfrak{a}f$, which are very ample (\cite{Har77} Ex.\ V.2.11). 
     Thus $L-K$ is very ample. 
   \end{proof}

   Next, let us describe which linear systems $L$ on a ruled surface
   $S$ determine a triple-point defective pair $(S,L)$. 

   We will show that example \ref{ex:ruled} 
   describes, in most cases, the only possibilities.
   In order to do so we first have to consider the possible
   algebraic classes of irreducible curves with 
   self-intersection zero on a ruled surface.

   \begin{lemma}\label{lem:ruled}
     Let $B\in|bC_0+b'f|_a$ be an irreducible curve
     with $B^2=0$ and $\dim|B|_a\geq 0$, then we are in one of the following cases:
     \begin{enumerate}
     \item[(a.1)] $B\sim_a f$,
     \item[(a.2)] $e=0$, $b\geq 1$, $B\sim_a bC_0$, and
       $|B|_a=|B|_l$, or
     \item[(a.3)] $e<0$, $b\geq 2$, $b'=\frac{b}{2}e<0$, $B\sim_a
       bC_0+\frac{b}{2}ef$ and $|B|_a=|B|_l$.
     \end{enumerate}
     Moreover, if $b=1$, then $S\cong C_0\times\pp^1$.
   \end{lemma}
   \begin{proof}
     See \cite{Kei01} App.\ Lemma G.2.
   \end{proof}     

   We can now classify the triple-point defective linear systems on a
   geometrically ruled surface. In order to do so we should recall the
   result of \cite{CM07} Prop.\ 18.

   \begin{proposition}\label{prop:length3}
     Suppose that, with the notation
     in \eqref{eq:alphabeta}, $\alpha$ is surjective,
     and suppose that $L$ and $L-K$ are
     very ample with $(L-K)^2>16$. Moreover, suppose that for $p\in S$
     general and for $L_p\in|L-3p|$ general the equimultiplicity
     scheme $Z_p$ of $L_p$ has a subscheme $Z_p'$ 
     of length $3$ such that  $h^1\big(S,\kj_{Z'_p}(L)\big)\not=0$.

     Then for $p\in S$ general there is an irreducible, smooth, rational curve $B_p$
     in a pencil $|B|_a$ with $B^2=0$, $(L-K).B=3$ and $L-K-B$ big.

     In particular, $S\rightarrow |B|_a$ is a ruled surface and $2B_p$ is a 
     fixed component of $|L-3p|$.     
   \end{proposition}

   \begin{theorem}\label{thm:ruled}
     With the above notation let $\pi:S\rightarrow C$ be a
     geometrically ruled surface, and let $L$ be a line bundle
     on $S$ such that $L$ and $L-K$ are very ample. Suppose that
     $(L-K)^2>16$ and that for a
     general $p\in S$ the linear system $|L-3p|$ contains a curve
     $L_p$ such that
     $h^1\big(S,\kj_{Z_p}(L)\big)\not=0$ where $Z_p$ is the
     equimultiplicity scheme of $L_p$ at $p$.

     Then $L=C_0+\mathfrak{b}\cdot f$ for some divisor $\mathfrak{b}$
     on $C$ such that $\mathfrak{b}+\mathfrak{e}$ is very ample and
     $|L-3p|$ contains a fibre of $\pi$  as fixed component with multiplicity two.     
     Moreover, if $e\geq -1$ then $\deg(\mathfrak{b})\geq 2g+1$ and we
     are in the situation of Example \ref{ex:ruled}.
   \end{theorem}
   \begin{proof}
     As in the proof of \cite{CM07} Thm.\ 19, since the case in which
     the length of $Z_p$ is $4$ has been ruled out in Remark \ref{remlength},
     we only have to consider the situations in  Proposition
     \ref{prop:length3} above. 

     Using the notation there we have a divisor $A:=L-K-B\sim_a aC_0+a'f$
     and a curve $B\sim_a bC_0+b'f$ satisfying certain numerical
     properties, in particular $p_a(B)=0$, $B^2=0$, and  $a>0$ since $A$ is big. Moreover, 
     \begin{equation}\label{eq:ruled:1}
       3=A.B=-eab+ab'+a'b
     \end{equation}
     and
     \begin{equation}\label{eq:ruled:2}
       a\cdot(2a'-ae)=A^2
       =(L-K)^2-2\cdot A.B-B^2
       \geq 17-2\cdot A.B-B^2=11.
     \end{equation}
     
     By Lemma \ref{lem:ruled} there are three possibilities for $B$ to
     consider.
     If $e<0$ and $B\sim_a bC_0+\frac{eb}{2}\cdot f$ with $b\geq 2$, then 
     Riemann-Roch leads to the impossible equation
     \begin{displaymath}
       -2=2p_a(B)-2=B.K
       =(2g-2)\cdot b.
     \end{displaymath}
     If $e=0$ and $B\sim_a bC_0$, then similarly Riemann-Roch shows
     \begin{displaymath}
       -2=B.K=(2g-2)\cdot b,
     \end{displaymath}
     which now implies that $b=1$ and $g=0$. In particular, $S\cong
     \pp^1\times\pp^1$ and $L\sim_a A+B+K\sim_a (a-1)\cdot C_0+f$,
     since $3=A.B=a'$. But this is then one of the cases of Example
     \ref{ex:ruled}. 

     Finally, if $B\sim_a f$ then \eqref{eq:ruled:1} gives $a=3$, and thus
     \begin{displaymath}
       L\sim_a A+B+K \sim_a C_0+(\mathfrak{a}'+\pi(p)+K_C+\mathfrak{e})\cdot f,
     \end{displaymath}
     where $A=3C_0+\mathfrak{a}'\cdot f$.
     Moreover, by the assumptions of Case (b) the linear system
     $|L-3p|$ contains the fibre of the ruling over $p$ as double
     fixed component, and since $L$ is very ample it induces on $C$
     the very ample divisor
     $\mathfrak{e}+(\mathfrak{a}'+\pi(p)+K_C+\mathfrak{e})$. 
     Note also, that \eqref{eq:ruled:2} implies that
     \begin{displaymath}
       a'-2-e\geq \frac{e}{2},
     \end{displaymath}
     and thus for $e\geq -1$ we have
     \begin{displaymath}
       \deg(\mathfrak{a}'+\pi(p)+K_C+\mathfrak{e})=2g+1+(a'-2-e)\geq 2g+1,
     \end{displaymath}
     so that then the assumptions of Example \ref{ex:ruled} are
     fulfilled. 
     This finishes the proof.
   \end{proof}

   If $\pi:S\longrightarrow C$ is a ruled surface, then there is a
   (not necessarily unique (if $g(C)=0$)) minimal model 
   
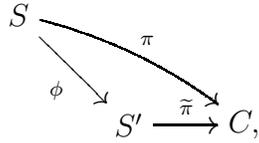
\begin{displaymath}
     \xymatrix{
       S\ar[dr]_{\phi}\ar@/^/[drr]^{\pi}& \\
       & S'\ar[r]^{\widetilde{\pi}} &C,
     }     
   \end{displaymath}
   and the N\'eron-Severi group of $S$ is\label{page:notation}
   \begin{displaymath}
     \NS(S)=C_0\cdot \Z\oplus f\cdot\Z\oplus\bigoplus_{i=1}^kE_i\cdot\Z,
   \end{displaymath}
   where $f$ is a general fibre of $\pi$, $C_0$ is the total
   transform of section of $\widetilde{\pi}$, and the $E_i$ are the
   total transforms of the exceptional divisors of the blow-up
   $\phi$. Moreover, for the Picard group of $S$ we just have to
   replace $f\cdot\Z$ by $\pi^*\Pic(C)$. We may, therefore,
   represent a divisor class $A$ on $S$ as
   \begin{equation}\label{eq:L}
     L=a\cdot C_0+\pi^*\mathfrak{b}-\sum_{i=1}^kc_iE_i.
   \end{equation}

   \begin{corollary}\label{cor:geometricallyruled}
     Suppose that $(S,L)$ is a tuple as in Proposition
     \ref{prop:length4} with ruling $\pi:S\rightarrow C$, and
     suppose that the N\'eron-Severi gruop of $S$ is as described before
     with general fibre $f=B_p$. 

     Then $S$ is minimal, 
     $L=C_0+\pi^*\mathfrak{b}$ for some divisor $\mathfrak{b}$
     on $C$ such that $\mathfrak{b}+\mathfrak{e}$ is very ample and
     $|L-3p|$ contains a fibre of $\pi$  as fixed component with multiplicity two.     
   \end{corollary}
   \begin{proof}
     Let $L=C_0+\pi^*\mathfrak{b}-\sum_{i=1}^k c_iE_i$, as described in \eqref{eq:L}. Then
     \begin{displaymath}
       L-K=(a+2)\cdot
       C_0+\pi^*(\mathfrak{b}-K_C-\mathfrak{e})-\sum_{i=1}^k(c_i+1)\cdot E_i,
     \end{displaymath}
     and thus considering Proposition \ref{prop:length3}
     \begin{displaymath}
       3=(L-K).B=a+2.
     \end{displaymath}
     
     The very ampleness of $L$ implies thus that $c_i>0$ for all $i$. But then, if
     $S$ is not minimal and $f'$ is the strict transform of a fiber of the minimal model,
     meeting some $E_i$, then $L\cdot f'\leq 0$, a contradiction.
   \end{proof}

   By \cite{CM07} we get Theorem
   \ref{main2} as an immediate corollary.


   \bibliographystyle{amsalpha-tom}

\begin{thebibliography}{BHPV04}

   \bibitem[BHPV04]{BHPV04}
     Wolf Barth, , Klaus Hulek, Christian Peters, and Antonius {Van de Ven},
     \emph{Compact complex surfaces}, 2nd ed., Ergebnisse der Mathematik und ihrer
     Grenzgebiete, 3.\ Folge, vol.~4, Springer, 2004.

   \bibitem[ChM07]{CM07}
     Luca Chiantini and Thomas Markwig, \emph{Triple-point defective regular
       surfaces}, arXive:0705.3912, May 2007.

   \bibitem[Fri98]{Fri98}
     Robert Friedman, \emph{Algebraic surfaces and holomorphic vector bundles},
     Springer, 1998.

   \bibitem[FuP00]{FuP00}
     Luis Fuentes and Manuel Pedreira, \emph{The projective theory of ruled
       surfaces}, math.AG/0006204, June 2000.

   \bibitem[Har77]{Har77}
     Robin Hartshorne, \emph{Algebraic geometry}, Springer, 1977.

   \bibitem[Kei01]{Kei01}
     Thomas Keilen, \emph{Families of curves with prescribed singularities}, Ph.D.
     thesis, Universit\"at Kaisers\-lautern, 2001, http:// \!\!www.
     \!\!mathematik. \!\!uni-kl. \!\!de/ \!\!\textasciitilde keilen/ \!\!download/
     \!\!Thesis/ \!\!thesis.ps.gz.

   \bibitem[LaB92]{LB92}
     Herbert Lange and Christina Birkenhake, \emph{Complex abelian varieties},
     Grundlehren der mathematischen Wissenschaften, no. 302, Springer, 1992.

   \bibitem[Rei97]{Rei97}
     Miles Reid, \emph{Chapters on algebraic surfaces}, Complex Algebraic Geometry
     (J\'anos Koll\'ar, ed.), IAS/Park City Mathematics Series, no.~3, Amer.\
     Math.\ Soc., 1997, pp.~3--159.

   \end{thebibliography}

   \providecommand{\bysame}{\leavevmode\hbox to3em{\hrulefill}\thinspace}

\end{document}